\documentstyle[harvard]{jsc}
\newtheorem{Algorithm}{Algorithm}[section]
 \title{   A new algorithm for computing the asymptotic solutions of a class of linear differential systems}
\author[B.M.Brown, M.S.P.Eastham, D.K.R.M$^c$Cormack]{B.M.Brown \and  M.S.P.Eastham \and D.K.R.M$^c$Cormack\\ \\
 Department of Computer Science, Cardiff University of Wales, Cardiff, CF2 3XF, U.K.}
\date{}
\begin{document}
\maketitle
 \newcommand{\beq}{\begin{equation}}
\newcommand{\enq}{\end{equation}}
\newcommand{\s}{ {\cal{x}}}
\large
\section{Introduction}
In this paper  we discuss a new algorithm for estimating and improving error terms in the asymptotic solution of linear differential systems. We consider systems of the form
\beq
Z^\prime(x) = \rho(x)  \{ D +R(x) \} Z(x) \;\;\; (a \leq x < \infty),\label{eq:1.1}
\enq
where $Z$ is an $n-$ component  vector, $ \rho$ is a real or complex scalar factor, $D$ is a constant $n \times n$ diagonal matrix,
\beq
D= dg(d_1,...,d_n) \label{eq:1.2}
\enq
with distinct $d_k$, and $R$ is also  an $n \times n$ matrix whose entries tend to zero as $ x \rightarrow \infty$, that is, $R(x) \rightarrow 0$ as $ x \rightarrow \infty$.
\par
If it is the case that $\rho (x) R(x)$ is $L(a, \infty)$, the Levinson asymptotic theorem \cite[section 1.3]{MSPE89};\cite{NL48}
states that there are solutions $Z_k$  $(1 \leq k \leq n)$ of (\ref{eq:1.1}) such that
\par
\beq
Z_k(x)=  \{e_k + \eta_k(x) \}\exp (d_k \int_a^x \rho( t)dt  ), \label{eq:1.3}
\enq
where $e_k$ is the unit coordinate vector in the $k-$ direction and   $\eta_k(x) \rightarrow 0 $ as $ x \rightarrow  \infty$. The size of the error term $\eta_k$ is related to the size of $R(x)$ as $x \rightarrow \infty$, and therefore the accuracy of ( \ref{eq:1.3})  can be improved if the perturbation $R(x)$ can be improved -- that is, made smaller in magnitude --- as $x \rightarrow \infty$. Under suitable conditions on
$\rho$ and $R$, this improvement can be effected by applying a sequence of transformations to the solution vector $Z$ in (\ref{eq:1.1}). Our algorithm is concerned with the implementation of this sequence of transformations.
\par
In order to introduce the ideas involved, we consider the transformation
\beq
Z =(I+P )W,    \label{eq:1.4}
\enq
where   $I$ is the   $ n \times n$ identity matrix, dg$P=0$,  and the off-diagonal entries of $P$ are defined by
\beq
PD-DP=R-dgR. \label{eq:1.5}
\enq
Thus, in terms of the $(i,j)$ entries in the matrices,
\beq
p_{ij} = r_{ij}/(d_j-d_i) \;\;\; ( i \neq j ). \label{eq:1.6}
\enq
On substituting (\ref{eq:1.4}) into ( \ref{eq:1.1}) and using ( \ref{eq:1.5}), we have
\beq
W^\prime = \rho \{ \tilde{D} +(I +P)^{-1} (RP-P {\rm dg} R - \rho ^{-1}P^\prime ) \} W, \label{eq:1.7}
\enq
where 
\beq \tilde{D} = D + {\rm dg} R. \label{eq:1.8}
\enq
By ( \ref{eq:1.6}), $ P(x) \rightarrow 0$ as $ x \rightarrow \infty $ and therefore it is clear that there are
circumstances, to be detailed later, in which the perturbation term in the $W$-system has a smaller order of magnitude for large $x$ than the original perturbation $R$. Repetition of the process successively improves the perturbation term.  It is to the final system in the process that the Levinson theorem ( \ref{eq:1.3}) is applied, when a prescribed accuracy in the error term has been achieved.
\par
The transformation back from the final system to the original system (\ref{eq:1.1}) yields an improvement of (\ref{eq:1.3})  in which the factor $e_k + \eta_k(x)$ is replaced by
\beq
\{ I+P_0(x) \}\{e_k + \eta_k(x)\} \label{eq:1.8a}
\enq
with a new $\eta_k$ which has the prescribed accuracy, and the matrix $P_0$ is generated explicitly by our algorithm. The terms in $P_0$ tend to zero as $x \rightarrow \infty$ and $ \eta_k = o(P_0)$. Thus (\ref{eq:1.8a})
provides explicit sub-dominant terms for the asymptotic solution of (\ref{eq:1.1}). In sections 2 and 3, we discuss the sequence of transformations and, in sections 4 and 5, we discuss the algorithm for the generation of the terms in $P_0$.
\par
 In a recent paper \cite{BEEM95}, we consider a particular system of the form ( \ref{eq:1.1}) which arises from the $n-$th order differential equation
\beq
y^{(n)} (x) - Q(x)y(x)=0, \label{eq:1.9}
\enq
and we formulated an algorithm which implements a sequence of transformations of the type (\ref{eq:1.4})--(\ref{eq:1.6}). The main emphasis in \cite{BEEM95}, however, is on the analytic and asymptotic implications of the  transformations for the solution of (\ref{eq:1.9}) and for applications to spectral theory. Here, on the other hand, we wish to develop our algorithm from the point of view of symbolic algebra in the context of the general system (\ref{eq:1.1}). We also demonstrate the versatility of our procedure by applying it to other situations than the one covered in \cite{BEEM95}.
\par
Finally in this introduction, we note that the origins of the transformation (\ref{eq:1.5})-(\ref{eq:1.6})
lie in the work of Harris and Lutz \cite{HL74} with more recent developments of these ideas by Eastham \cite{MSPE86}\cite[section 1.7]{MSPE89}. The nature of the matrix $I+P$ is that it is an explicit approximation to the matrix whose columns are eigenvectors of $D+R$, these eigenvectors in general only being explicit in terms of $D$ and $R$ when $n=2$. 

\section{The sequence of transformations}
We now define a sequence of transformations
\beq
Z_m =(I+P_m)Z_{m+1}   \;\;\;(m=1,2,...) \label{eq:2.1}
\enq
of the type introduced in (\ref{eq:1.4}) -- (\ref{eq:1.6}), the purpose of
which is to produce differential systems for the $Z_m$, similar to 
(\ref{eq:1.1}), but with the perturbation term successively improved. The definition
is almost the same as that given in \cite[section 4]{BEEM95} for the particular 
system (\ref{eq:1.1}) which arises from (\ref{eq:1.9}), and so we shall
be brief in this part of the paper. A typical system in the process is
\beq
Z_m^\prime = \rho (D_m + R_m)Z_m, \label{eq:2.2}
\enq
where $D_m$ is diagonal, with (\ref{eq:1.1}) being the case $m=1$.
The process ends at $m=M$ when the perturbation $R_M$ has a pre-assigned
accuracy in terms of its order of magnitude as $x \rightarrow \infty$.
\par
As already indicated by (\ref{eq:1.7}), $R_m$ will contain  
terms of different  orders of magnitude as $x \rightarrow \infty$, 
and it is the essence of our algorithm to identify and collate 
these terms according to their size. Hence we write
\beq
R_m = V_m +E_m = V_{1m} + V_{2m} +...+V_{\mu m}+E_m, \label{eq:2.3}
\enq
where
\beq
V_{km}=o(V_{jm}) \;\;\;( x \rightarrow \infty, k >j)  \label{eq:2.4}
\enq
and
\beq
E_{m}=o(V_{\mu m}) \;\;\;(x \rightarrow \infty).  \label{eq:2.5}
\enq
Here $E_m $ represents terms which are already of the pre-assigned
accuracy, and they take little part in the transformation process
(\ref{eq:2.1}). The $V_{jm}$ represents terms which are not of that
accuracy, and they are successively replaced by smaller-order 
terms as we go through the process.  Also as indicated by (\ref{eq:1.8}),
we take any diagonal terms in $V_{1m}$ over to $D_m$ in (\ref{eq:2.2}). 
Thus we arrange that
\beq
{\rm dg} V_{1m}=0 \label{eq:2.6}
\enq
and we write
\beq
D_m =D+\Delta_m.
\enq
\par
We substitute (\ref{eq:2.1}) into (\ref{eq:2.2}) to eliminate the 
dominant term    $V_{1m}$ in (\ref{eq:2.3}) and to define the 
resulting terms $V_{j,m+1}$ constructively in terms of the $V_{jm}$. 
Corresponding   to (\ref{eq:1.5}), we define $P_m$ by
\beq
P_mD-DP_m = V_{1m} \label{eq:2.8}
\enq
with dg$P_m=0$. Then it is easily checked that (\ref{eq:2.1}) and
(\ref{eq:2.2}) give
\begin{eqnarray}
Z^\prime_{m+1} = \rho \{ D_m +(I+P_m)^{-1}(&-&\rho^{-1} 
P^\prime_m  + T_m +V_{1m}P_m  \nonumber \\
&+&(R_m-V_{1m})(I+P_m)) \}Z_{m+1},            \label{eq:2.9}
\end{eqnarray}
where
\beq
T_m=\Delta_m P_m -P_m \Delta_m.      \label{eq:2.10}
\enq
As in \cite[ section 4]{BEEM95}, we show that (\ref{eq:2.9}) can be expressed as
\beq
Z^\prime_{m+1} = \rho ( D_{m+1} + R_{m+1} ) Z_{m+1} \label{eq:2.11}
\enq
where $R_{m+1}$ has the form
\beq
R_{m+1}=V_{m+1}+E_{m+1}=V_{1,m+1}+...+V_{\mu, m+1}+E_{m+1} \label{eq:2.12}
\enq
as in (\ref{eq:2.3})-(\ref{eq:2.5}), but with a different $\mu$. To do this,
we let $U$ denote any of the terms on which $(I+P_m)^{-1}$ acts 
in (\ref{eq:2.9}). Then we write
\beq
(I+P_m)^{-1} = I-P_m+P^2_m- ... +(-1)^\nu P^\nu_m + (-1)^{\nu+1}(I+P_m)^{-1}
P^{\nu+1}_m   \label{eq:2.13}
\enq
where, for each $U, \nu$ is chosen so that the product
\beq
(I+P_m)^{-1}P^{\nu +1}_m U  \label{eq:2.13a}
\enq
has a sufficiently small order of magnitude to be included with $E_m$
and form part of $E_{m+1}$. Now we group together terms of the same order
of magnitude and denote the dominant term by $S_{m+1}$. 
We then obtain (\ref{eq:2.12}) (with $S_{m+1}$ in place of $V_{1,m+1}$),
where $E_{m+1}$ has the pre-assigned accuracy and, by (\ref{eq:2.8}),
$S_{m+1}$ and the $V_{j,m+1}$ are known explicitly in terms of the $V_{jm}$.
Then, finally, we obtain (\ref{eq:2.11}) and (\ref{eq:2.12}) by defining
\begin{displaymath}
D_{m+1}=D_m + {\rm dg}S_{m+1}
\end{displaymath}
and 
\begin{displaymath}
V_{1,m+1}=S_{m+1} - {\rm dg}S_{m+1}.
\end{displaymath}

 \section{Orders of magnitude}
 The transformation process (\ref{eq:2.1}) is carried out for $m=1,2,...,M-1$ 
 and, in order    to express the process in terms of an algorithm which
 can be implemented in the symbolic algebra system Mathematica, it is
 necessary to specify more precisely the orders of magnitude involved.
 The starting point is $m=1$ and, in (\ref{eq:2.3}), we suppose that
 \beq
 V_1=V_{11} + V_{21} +...+V_{N1}, \label{eq:3.1}
 \enq
 where
 \beq
 V_{j1}(x) =O(x^{-\theta_j}) \;\;\;(1 \leq j\leq N) \label{eq:3.2}
 \enq
 as $x \rightarrow \infty$ and, corresponding to (\ref{eq:2.4}),
 \begin{displaymath}
 0< \theta_1 < \theta_2 <...< \theta_N.
 \end{displaymath}
 We assume that the $\theta_j$ in (\ref{eq:3.2}) are chosen to have their
 minimum possible values and, in practice, (\ref{eq:3.2}) 
 represents the exact order of magnitude of $V_{j1}$. We denote by $\sigma$
 the set of positive numbers
 \beq
 \sigma = \{ n_1 \theta_1 + n_2 \theta_2+...+n_N \theta_N; \; n_1 \geq 1, \;
 n_2 \geq 0, ..., n_N \geq 0 \}     \label{eq:3.3}
 \enq
 the $n_j$ being integers. It is possible that different values of 
 the $n_j$ give   the same number in $\sigma$ and, allowing for this, we denote
 the distinct numbers in $\sigma$ by $\sigma_1,\sigma_2,...$ in 
 increasing order. Let us suppose that the pre-assigned accuracy represented 
 by $E_m$ in (\ref{eq:2.3}) is expressed as
 \beq
 E_m(x)=O(x^{-K})  \label{eq:3.4}
 \enq
 for a given $K>0$. Then we choose the integer $L$ so that
 \beq
 \sigma_L < K  \leq \sigma_{L+1}.  \label{eq:3.5}
 \enq
 \par
 The definition of $\sigma$ in (\ref{eq:3.3}) allows us to 
 postulate orders of magnitude
 \beq
 V_{jm}(x)=O(x^{-\sigma_{m+j-1}}) \label{eq:3.6}
 \enq
 where we can allow the possibility that some of the $V_{jm} $ 
 ( even $V_{1m} $) may be identically zero.   To justify (\ref{eq:3.6}),
 we note first that
 $P_m=O(x^{-\sigma_m})$ by (\ref{eq:2.8}). Then, recalling the use of 
 (\ref{eq:2.13}) in ( \ref{eq:2.9}), we also have
 \begin{displaymath}
 P^r_m V_{jm} = O(x^{-r \sigma_m - \sigma_{m+j-1}}),
 \end{displaymath}
 and again $ r \sigma_m + \sigma_{m+j-1} \in \sigma$ by (\ref{eq:3.3}).
 Further, since the combination $r=0$ and $j=1$ does not occur together
 here, we have
 \begin{displaymath}
 r \sigma_m + \sigma_{m+j-1} \geq \sigma_{m+1}.
 \end{displaymath}
 The term $T_m$ in (\ref{eq:2.9}) is treated similarly. A simple induction 
 argument on $m$ now establishes (\ref{eq:3.6}) for all $j$ and $m$,
 provided only that we add a suitable hypothesis on the term 
 $\rho^{-1}P^\prime_m$ which appears in (\ref{eq:2.9}) but is not so far included 
 in the argument.     We therefore add the hypothesis that
 \beq
 \rho^{-1}P^\prime_m = W_{1m} +...+W_{lm}      \label{eq:3.7}
 \enq
 where, similarly to (\ref{eq:3.6}),
 \beq
 W_{jm}(x)=O(x^{-\sigma_{m+j}}) \label{eq:3.8}
 \enq
 and again  we allow the possibility that some $W_{jm}$ may be zero. Since
 $P_m^\prime$ depends on $V^\prime_{1m}$ ( see (\ref{eq:2.8})), which 
 in turn depends   on the previous matrices in the process (\ref{eq:2.1}),
  the nature of (\ref{eq:3.7}) and ( \ref{eq:3.8}) is that they are
  conditions on the successive derivatives of the original $V_{j1}$ which
  occur in $R_1$ in (\ref{eq:2.3}) and (\ref{eq:3.1}). 
The exact form of these conditions on $V_{j1}$ determines classes of matrices $R_1$ to which
this theory and the consequent algorithms are applicable.  Examples of such classes will be given in
section 5.
Thus (\ref{eq:3.7})
  and (\ref{eq:3.8}) are consequences of these conditions on $R_1$   which
  must be established ( usually by induction) in each application of
  the theory. It is these $W_{jm}$ which will appear in   our algorithms.
  \par
  We can now summarise this section by saying that, subject to (\ref{eq:3.2}),
  (\ref{eq:3.7}) and ( \ref{eq:3.8}), we have established that
  \begin{displaymath}
  V_{jm}(x) = O(x^{-\sigma_{m+j-1}})
  \end{displaymath}
  in (\ref{eq:2.3})-(\ref{eq:2.5}). Also, allowing for the fact that some
  $V_{jm}$ in (\ref{eq:2.5}) may be zero, we can write 
  $\mu=L-m+1$ by (\ref{eq:3.4})  and (\ref{eq:3.5}). The transformation
  process (\ref{eq:2.1}) ends when (\ref{eq:2.3}) reduces to
  \beq
  R_M =E_M=O(x^{-K}), \label{eq:3.9}
  \enq  the pre-assigned accuracy, and it follows from (\ref{eq:3.5}) that
  \begin{displaymath}
  M=L+1, \mu = M-m.
  \end{displaymath}

 \section{The algorithm }
In this section of the paper we show how the theory that has been developed in sections 1 through 3  may be used to obtain a computer code to calculate the asymptotic expansion of the solutions of (\ref{eq:1.1}) together with an explicit  error bound at some point $x \geq X > 0$.
A consequence of the theory that we have exhibited is that, given sufficient computational power,  
the quality of the asymptotics  that we obtain allows us to take $X$
to be quite small and still maintain high accuracy in the solutions.
 \par
As in the discussion in \cite{BEEM95} the    algorithm is formulated and 
implemented in three stages.   All the symbolic algorithms that we shall discuss are implemented in the
symbolic algebra system Mathematica.
The first algorithm, which is concerned with the generation of  a set of
recurrence relations to compute the matrix quantities $S_j$,
assumes only that the quantities involved satisfy 
non-commutative multiplication.
We recall the comments made after (\ref{eq:3.8}) that general classes of matricies
$R_1$ to which the algorithm is applicable will be given in section 5.
 In the following, we  write $A_m=(I+P_m)^{-1}$  and we note that expressions such as $P_m$, $T_m$ and $W_{jm}$ appear in the algorithm by virtue of their orders
of magnitude as indicated in section 3.

\begin{Algorithm}

\newcounter{rem1}
\newcounter{rem2}
\newcounter{rem3}
\begin{list}%
{( \alph{rem1} )}{\usecounter{rem1}
\setlength{\rightmargin}{\leftmargin}
\setlength{\rightmargin}{\labelwidth}
\setlength{\leftmargin}{\labelwidth}}
\item
Define $K$ to specify the accuracy (\ref{eq:3.4}).
\item
Define $N$ and $\theta_1, ..., \theta_N$ in (\ref{eq:3.2}) and arrange the distinct numbers in the set $\sigma$ in increasing order. This defines 
$\sigma_1, \sigma_2,...$ and determines $L$ in (\ref{eq:3.5}).
For a given $K$, $n_j$ in (\ref{eq:3.3}) satisfies
\begin{list}%
{( \Roman{rem2} )}{\usecounter{rem2}}
\item
$0 \leq n_j \leq [ \frac{ K-\theta_1}{\theta_j}] \;\; (j \geq 2 )$
\item
$1 \leq n_1 \leq [ \frac{ K }{\theta_1}]$.
\end{list}
\item
Start with $D_1$ and $V_{j1} \; ( 1 \leq j \leq N) $ as in (\ref{eq:3.1})--(\ref{eq:3.2}) and put $E_1=0$.
\item
For   $ m=1$ to $M-1$, 
\begin{list}%
{( \Roman{rem2} )}{\usecounter{rem2}}
\item
$E_{m+1}=A_mE_m(I+P_m)$.
\item
For each $U \in \{W_{jm} \; ( 1\leq j \leq l), T_m,V_{1m}P_m,V_{jm}, V_{jm}P_m \;\; (2 \leq j \leq M-m )\}$, 
\begin{list}%
{( \roman{rem3} )}{\usecounter{rem3}}
\item
In (\ref{eq:2.13a}) determine $\nu$.
\item
For $r=0$ to $\nu$, 
\begin{itemize}
\item
determine the order $\sigma_{m+k} =r\sigma_m+($order of $U)$ of $P_m^r U$;
\item
Update $V_{k,m+1}= V_{k,m+1} + (-1)^rP^r_m U$.
\end{itemize}
\item
Update $E_{m+1}=E_{m+1}+(-1)^\nu A_mP_m^{\nu+1} U$.  
\end{list}
\item
Output $S_{m+1} = V_{1,m+1}$.
\end{list}
\end{list}
\end{Algorithm}
\par
At each stage of the algorithm, $S_{m+1}$ depends on the terms $D_k,P_k,W_{jk}$ and $V_k$ $(1 \leq k \leq m )$. However, because of part $(b)$, the algorithm requires more precise information than its counterpart, Algorithm 6.1, in \cite{BEEM95}. The set $\sigma$ in \cite{BEEM95} has a very simple form,  consisting only of numbers $na$, where $n$ is a positive integer and $a \; (>0)$ is a parameter. Thus 
$\sigma_n = na$ in (\ref{eq:3.3}), and Algorithm 6.1 in \cite{BEEM95} can be executed without specifying the value of $a$. We give a more general example of the same situation in Example 5.1 below. However, in the wider context of (\ref{eq:3.3}), sufficient information about the parameters $\theta_1,...,\theta_N$ must be provided to Algorithm 4.1 in order to generate all the necessary values of $\sigma_n$. We therefore defer further discussion of the output of Algorithm 4.1 to Example 5.2 in the next section, where values of the parameters are specified.
\par
  \begin{Algorithm}
Starting with  the precise form of the matrices $D_1$ and $V_1$ and with
 $E_1=0$, the expressions $S_2,...,S_{M-1}$ generated by Algorithm 4.1 
are evaluated in order. These are then used to evaluate the matrices 
$D_{m+1}$ and $V_{1,m+1}$.
\end{Algorithm}
  The structure of the algorithm is similar to that of Algorithm 6.2 of \cite{BEEM95}.
As noted in that paper, in order to reduce the computation time a detailed assessment
of the mathematical issues involved at each simplification of an expression must first be made.
Thus the judicious use of the {\it Together, Apart} commands instead of the {\it Simplify}
command can result in a dramatic  decrease in the time needed to perform the computation.
At this stage it is necessary to keep the expressions in symbolic form since the $W_{jm}$
must be obtained explicitly. These are computed in terms of $P^{'}_m$, which in turn is obtained from $S_j$ by differentiation.
\par
The final algorithm that is needed in the symbolic part of the computation 
is concerned with obtaining an upper bound for the norm $\parallel E_m \parallel
$   of the error term $E_m$. The norm is computed using the sup. norm
by applying the    triangle and Cauchy inequalities
\begin{displaymath}
\parallel AB \parallel \leq n \parallel A \parallel\parallel B \parallel , \;
 \parallel A+B \parallel \leq \parallel A \parallel+\parallel B \parallel
 \end{displaymath}
 for $ n \times n $ matrices.   As in \cite{BEEM95}  
 a bound for the norm of the inverse matrix $A_m= (I+P_m)^{-1}$ is given by
 \beq
 \parallel A_m \parallel \leq 1 + \parallel P_m \parallel/(1-n
 \parallel P_m \parallel )       \label{eq:4.4a}
 \enq
 provided $ \parallel P_m \parallel <1/ n $.
\begin{Algorithm} \label{alg:6.3}
Compute the ${\it sup.}$ norm of each matrix in $E_m$, using (\ref{eq:4.4a})
for the inverse matrices.
Next apply the triangle and Cauchy inequalities to obtain an upper bound for 
the ${\it sup.}$ norm of $E_m$ itself. 
 \end{Algorithm}
We note that, since Algorithm 4.1 expresses $E_m$ in terms of matrices arising at earlier stages of 
(\ref{eq:2.2})-- and therefore ultimately in terms of the $V_{j1}$ in(\ref{eq:1.1})--
so also Algorithm 4.3 ultimately expresses the norm of $E_m$ in terms of norms derived from the original system (\ref{eq:1.1}).
\par
 Before moving on to examples of the implementation of the algorithms, we add some detail to (\ref{eq:1.8a}) concerning the generation of the sub-dominant terms in the asymptotic solution of (\ref{eq:1.1}). By (\ref{eq:3.9}) and (\ref{eq:3.5}), the final system in the sequence (\ref{eq:2.2}) is
\beq
Z'_M = \rho (D_M+E_M)Z_M, \label{eq:4.2}
\enq
where
\beq
\parallel E_M \parallel \leq c_Mx^{- \sigma_M} \label{eq:4.3}
\enq
and $c_M$ is a constant. We recall that the numbers $\sigma_m$ cover all orders of magnitude which occur. Algorithm 4.3 provides a definite value for $c_M$ in any particular example. As in ( \ref{eq:1.3}), the asymptotic solution of (\ref{eq:4.2}) has the form
\par
\beq
\{ e_k + \eta_k(x) \} \exp ( \int_a^x d_{kM}(t) \rho(t) dt ), \label{eq:4.4}
\enq
where the $d_{kM}$ are the diagonal entries in $D_M$ and the size of $\eta_k$ can be expressed in terms of $c_M$ and $\sigma_M$ as in\cite [(3.15)]{BEEM95}. What we wish to emphasise here is the transformation back from (\ref{eq:4.2}) to the original system (\ref{eq:1.1}).
As indicated in (\ref{eq:1.8a}), this adds to (\ref{eq:4.4}) the extra factor
\beq
I+P_0(x) = \Pi_{m=1}^{M-1} \{ I + P_m(x)\}.
\label{eq:4.5}
\enq
Now the definition of $P_m$ in (\ref{eq:2.8}) is in terms of $V_{1m}$, which is provided by Algorithms 4.1 and 4.2. Further, by (\ref{eq:2.8}) and (\ref{eq:3.6}), $P_m(x) = O(x^{-\sigma_m}) \; ( 1 \leq m \leq M-1 ).$
Thus, in terms of (\ref{eq:4.5}), our algorithms provide sub-dominant terms up to 
$O(x^{-\sigma_{M-1}})$ 
in the asymptotic solution of (\ref{eq:1.1}).
\par
We mention one further point concerning the transformation process which leads from (\ref{eq:2.2}) to (\ref{eq:2.9}). Since the derivative $P'_m$ appears in (\ref{eq:2.9}) and since $P_m$ ultimately depends on $V_1$ and $\rho$, each step in the process requires the existence of a further derivative of $V_1$ and $\rho$. Thus the sub-dominant terms in (\ref{eq:4.5}) require  the existence of $M-1$ derivatives of $V_1$ and $\rho$. If $V_1$ and $\rho$ are infinitely differentiable then, subject to convergence considerations, (\ref{eq:4.5}) would yield a full asymptotic expansion. It is hoped to deal with this matter in a future paper.

  \section{Examples}
 \subsection{Example 1}
 Let $\rho(x)=x^\gamma$ and $ R(x)=x^{-(1+\gamma)}C$, where $\gamma > 0$  and $C$ is a constant matrix.
 Here we have just $N=1$ in  (\ref{eq:3.1}) and
\begin{displaymath}
\sigma_m=m(1+\gamma) \;\; ( m=1,2,...).
\end{displaymath}
This example is basically the case considered in \cite [section 5]{BEEM95} with a special choice of $C$ and, as in \cite{BEEM95}, the condition  (\ref{eq:3.7})
is easily verified by induction on $m$. The present code has been  tested on this example 
and the results  from Algorithm 4.1 are, up to notational differences, identical with those reported on in \cite{BEEM95}. Further, Algorithms 4.2 and 4.3   return  values of the solutions
 computed with $4$ iterations that, at $X=40$, are within $10^{-11}$ of those reported on in \cite{BEEM95}.
   \subsection{Example 2}
A significantly different example is obtained when $\rho(x)$ and  $R(x)$
in (\ref{eq:1.1}) contain periodic factors.
Let  $\rho(x)=x^\gamma  \phi(x^\beta)$
and 
\beq
R(x)=x^{-(1+\gamma-\beta)}F_1(x^\beta)   + x^{-(1+\gamma)}F_2(x^\beta),
\label{eq:5.1}
\enq
where
\beq
0< \beta < 1+\gamma \label{eq:5.2}
\enq
and $\phi(t),F_1(t)$ and $F_2(t)$ have the same period $\omega$ in $t$, with $\phi$ nowhere zero.
Here we have $N=2$ in (\ref{eq:3.1}) and
\beq
\theta_1=1+\gamma-\beta,\;\;\; \theta_2=1+\gamma \label{eq:5.3}
\enq 
in  (\ref{eq:3.1})-(\ref{eq:3.2}).
Corresponding to (\ref{eq:3.6}), we make the induction hypothesis
\begin{displaymath}
V_{jm}=x^{- \sigma_{m+j-1}}U_{jm}(x^\beta)
\end{displaymath}
where $U_{jm}(t)$ has period $\omega$ in $t$. Then, by  (\ref{eq:2.8}),
\begin{displaymath}
P_m(x)=x^{-\sigma_m} \Pi_m ( x^\beta)
\end{displaymath}
where $\Pi_m(t)$ has period $\omega$ and the entries $\pi_{ijm}$ in $\Pi_m$ are obtained from those in   $U_{1m}$ by the formula
\begin{displaymath}
\pi_{ijm}=u_{ij1m}/(d_j-d_i) \;\;\;( i \neq j ).
\end{displaymath}
Then  considering (\ref{eq:3.7}), we have
\begin{eqnarray}
\rho^{-1}P_m^{'} &=& x^{ -( \sigma_m + \gamma +1 -\beta)}(\Pi^{'}_m/\phi)(x^\beta) 
-\sigma_m x^{ -( \sigma_m + \gamma +1 )}(\Pi_m/\phi)(x^\beta) \nonumber \\
&=& x^{ -( \sigma_m + \gamma +1 -\beta)} W_1( x^\beta) + x^{ -( \sigma_m + \gamma +1 )}W_2(x^\beta) \label{eq:5.4}
\end{eqnarray} 
and hence (\ref{eq:3.7}) holds with $l=2$.
\par
We note that the upper bound (\ref{eq:5.2}) placed on $\beta$ is a restriction on the frequency of oscillations of $\rho$ and $R$ in this example. The  same type of condition was imposed in \cite[Example 2.4.1]{MSPE89}  in connection with the method of repeated diagonalization. When $\beta > 1 + \gamma$,  the asymptotic solution of (\ref{eq:1.1}) requires transformations of an entirely different nature from those based on (\ref{eq:2.1}) and (\ref{eq:2.8}) \cite{MSPE92a}, \cite{YTS}.
\par
In order to discuss the output of Algorithm 4.1 for this example, we have to choose specific values of $\beta$ and $\gamma$, so that part (b) can generate the list of values $\sigma_m$. We make the simple choice $\beta=\gamma=1$, so that  $\theta_1=1$ and $\theta_2=2$ in (\ref{eq:5.3}). Then, by (\ref{eq:3.3}), $ \sigma_m =m$.  
Also, by (\ref{eq:5.1}) and (\ref{eq:5.4}), we have
\begin{eqnarray*}
V_{11}(x) = x^{-1}F_1(x), && V_{21}(x) = x^{-2}F_2(x),  \\
W_{1m}(x)= x^{-(m+1)}W_1(x), && W_{2m}(x) = x^{-(m+2)}W_2(x)
\end{eqnarray*}
in the notation of section 3.
\par
 We have noted in section 4 that $S_{m+1}$ depends on $D_k, \; P_k, W_{jk}$ and $V_k \; (1 \leq k \leq m).$  However, the formulae for the $S_{m+1}$
 can often be simplified by expressing them in terms of the earlier $S_k$. This reduces the number of terms in the formulae with a consequent reduction in the computational effort required. We now give the output for $S_2,S_3$ and $S_4$:
\begin{eqnarray*}
S_2 &=& V_{11}  P_1  + T_1  + V_{12}  - W_{11}  \\ \nonumber
S_3 &=& -P_1 S_2 + V_{12}  P_1  + T_2  - W_{12}  - W_{21}  \\ \nonumber
S_4 &=& T_3  - W_{13} \\ \nonumber
\end{eqnarray*}
We note that, at this stage, these expressions appear no more complex that those computed in \cite{BEEM95}. However, increased difficulties do occur in the evaluation of the  entries in $T_3$, $W_{12},W_{21}$  and $W_{31}$ at  the next stage  when Algorithm 4.2 is implemented. We discuss this point further in Example 5.3. The time needed on a Sun SPARC-station 10 to compute $S_4$ is $1.6$ seconds, which is comparable with the comparative time $1.1833$ seconds reported on in \cite{BEEM95}. 
The similar times are a reflection of the low number of terms that must be manipulated by the symbolic algebra system.
As we remarked previously, the output of Algorithm 4.1 at this point consists only of a set of symbolic expressions which satisfy non-commutative multiplication. 
\par
The error term
\begin{eqnarray*}
E_4 &=& -A_3   W_{23}  
 + A_4  T_4  -  A_4   W_{14}    A_4
 W_{24}   \\ 
&+& A_1   P_1 ^2  \left( V_{12}  P_{1}  - W_{21}  \right) \\
&+& A_2  \left( -P_1  \left( T_1  + V_{12}  - W_{11}\right)+ V_{12}  P_1  - P_1  V_{11}  P_1  - W_{21}  \right)  P_2  \\
&+& A_3 V_{31}  P_3  + A_4  V_{41}  P_4 
\end{eqnarray*}
is more complex than that found in \cite{BEEM95}, which consists of only   $11$ additive terms. This additional complexity in the $E_4$ term is reflected in the time needed to compute norms when the third stage, Algorithm 4.3, is implemented.
\subsection{Example 3}
Algorithms 4.2 and 4.3 require the input of specific matrices $D_1$ and $V_1$ and, in this example, we introduce a special case of Example 5.2 which arises from the $n-$th order differential equation
\beq
y^{(n)}(x)-Q(x)y(x)=0. \label{eq:5.5}
\enq
Again, $Q(x)$ contains a periodic factor of the form 
 \beq
 Q(x)= x ^\alpha f(x^\beta)  \label{eq:5.6}
 \enq 
 where $f(t)$ is periodic in $t$ and nowhere
zero,  with $0 < \beta < 1 +   \frac{\alpha}{n}$.
As in \cite[section 3]{BEEM95}, we write (\ref{eq:5.5}) in the system form 
\beq
Z'=Q^{1/n}(D+Q'Q^{-1-1/n}C)Z, \label{eq:5.7}
\enq
where $Z$ has a first component $y$, $D$ is the diagonal matrix formed by
the $n-$th roots of unity,
\beq
D=\rm{dg}( \omega_1,...,\omega_n),
\enq
 and $C$ is constant with 
\beq
{\rm dg} C = -(n-1)(2n)^{-1}I. \label{eq:5.8}
\enq
It follows from (\ref{eq:5.6}) that
\begin{displaymath}
Q'Q^{-1-1/n}= \beta x^{-(1+\alpha/n-\beta)}( f'f^{-1-1/n})(x^\beta)
+ \alpha x^{-(1+\alpha/n)}f^{-1/n}(x^\beta).
\end{displaymath}
Hence (\ref{eq:5.7}) is the special case of (\ref{eq:1.1}) and (\ref{eq:5.1}) in which 
\begin{displaymath}
\gamma = \alpha/n,\; \phi=f^{1/n},\; F_1=\beta f'f^{-1-1/n}C,\; F_2 = \alpha f^{-1/n}C.
\end{displaymath}
\par
We now write (\ref{eq:5.7}) in the form (\ref{eq:2.2}) ( with $m=1$), where $dg V_{11}=0$ as in (\ref{eq:2.6}). Thus taking the diagonal terms from $F_1$ over to $D$ and
using (\ref{eq:5.8}), we define
\beq
D_1=D-(n-1)(2n)^{-1} \beta x^{-(1+\alpha/n-\beta)}(f'f^{-1-1/n})(x^\beta) I = D+\frac{1}{2}(n-1)pI, \label{eq:5.9}
\enq
where
\begin{displaymath}
p=x^{-\alpha/n} \{ f^{-1/n}(x^\beta) \}'
\end{displaymath}
and
\beq
R_1=V_{11}+V_{21} \;\;(=V_1), \label{eq:5.10}
\enq
where
\begin{displaymath}
V_{11}=-x^{-\alpha/n}np(C-{\rm dg } C), \; V_{21}=\alpha x ^{-(1+\alpha/n)}f^{-1/n}(x^\beta)C.
\end{displaymath}
Thus (\ref{eq:5.9}) and (\ref{eq:5.10}) are our choice of $D_1$ and $V_1$.
\par
As in the discussion of Example 5.2, we choose the parameter values $\beta=\gamma=1$, that is, $\alpha=n$ and $\beta=1$ in (\ref{eq:5.6}).
Finally, we must also choose the values of $n$  in order to complete the requirements for implementing Algorithms 4.2 and 4.3. 
We choose $n=4$, in which case a short calculation gives
\begin{displaymath}
C=-\frac{1}{8}\left ( \begin{array}{cccc}
-3 & 1+i & 1 & 1-i \\
1-i & -3  & 1+i & 1 \\
1 & 1-i & -3 & 1+i \\
1+i & 1 & 1-i & -3
 \end{array}
\right )
\end{displaymath}
as in \cite[Algorithm 6.2]{BEEM95}.
\par
The periodic nature of the function $Q(x)$ introduces 
additional matrices    over the case discussed in Example 5.1 and
the theory  expounded in \cite{BEEM95}.  As we remarked above, 
Algorithm 4.1 needs the specific values of $\beta$ and $\gamma$
 to be available. A   consequence of this is that all three algorithms must be run for each set  of parameter values. However the main additional computational difficulties occur in Algorithms 4.2 and 4.3. In Algorithm 4.2   the extra matrices generated as a consequence of the periodic nature of $Q$ must have their entries evaluated, while in Algorithm 4.3 the norm of the error matrix, which is considerably more complex than that which occurs in Example 5.1,
must be evaluated.  
\par
In order to test the performance of the set of algorithms  we have chosen 
to take
\beq
f(x)=2+\sin x. \label{eq:5.13}
\enq
However the performance of Algorithms 4.2 and 4.3 is considerably improved if we 
work with a generic function $f$ together with the simplification rule
\beq
f^{''}=2-f
\enq
and the results that we report are based on this latter situation.
A further consequence of the extra complexity in $Q$ is that, with the CPU power and memory that we have available, we can not evaluate the entries in $S_5$. Thus   effectively, we can only perform $4$ iterations of Algorithm 1.
  The time needed on a SPARC 10 workstation to compute the entries in $S_4$ is 
$190$ seconds of CPU time. This compares  with the $75$ CPU seconds that was needed in the work reported on in   \cite{BEEM95}.
The final algorithm, Algorithm 4.3, deals with the estimation of the $sup.$ norm of the error matrix, in our case $E_4$.
This involves first applying the Cauchy and triangle inequalities to each
matrix component of $E_4$ to obtain an upper bound for 
$ \parallel E_4 \parallel $ in terms of the norms of its components. An estimate for the norm of the inverse matrix $A_4$ is given by (\ref{eq:4.4a}).
\par
 These  norms   are estimated by examining  and evaluating each component of each matrix.
In doing this we encounter terms involving $p$ and its derivatives. 
In order to obtain upper bounds for these terms we   symbolically compute expressions  for them and note that, for this example, 
  \begin{displaymath}
3   \geq |f(x)| \geq   1  \;\; (X<x<\infty)    
 \end{displaymath}
This gives the   necessary bound.
 Again the increase in the complexity  of the expressions, resulting from 
  the more complex structure of the initial data means that the computation
  time that is required is increased over \cite{BEEM95}. 
  It takes some $1110$ CPU seconds to compute the norm of $E_4$ compared with    approximately $550$ CPU seconds  for $E_4$ reported in \cite{BEEM95}. 
At $X=40$ the bound for  the norm of this error matrix is $ 1.63099 \times 10^{-6}$.
  \par
We now compute the factor 
\begin{displaymath}
I+P_0(x)
\end{displaymath}
 discussed in (\ref{eq:4.5}) and apply this to the asymptotic solution (\ref{eq:4.4}) to yield the asymptotic solution of (\ref{eq:5.7}).
We mention that in view of the size of the expressions that are generated in 
(\ref{eq:4.5}), we have chosen to  evaluate it at $x=40$ using $30$ digits of accuracy.

\section{Concluding remarks}
\subsection{The transformations of Harris and Lutz}
In the introduction, we indicated that the origins of our basic transformation (\ref{eq:1.4})-(\ref{eq:1.6}) lie in the 1974 paper of Harris and Lutz.
In a subsequent paper (\cite[2.4]{HL77}), an extension of (\ref{eq:1.5}) is also discussed. Whereas (\ref{eq:1.5}) can be described as providing a first-order approximation to the exact diagonalization of $D+R$
in (\ref{eq:1.1}), the extension provides a more accurate second-order approximation. These ideas are also discussed in (\cite{MSPE89} pp26-8).
\par
The question arises whether this extension accelerates the process leading to (\ref{eq:3.9}),
and here we indicate why it does not achieve this objective. The essential feature of both (\ref{eq:1.5})
and the extension in \cite[2.4]{HL77} is that they are linear algebraic equations to determine $P$.
They do not  involve $P'$.
Thus, in (\ref{eq:2.9}), both the corresponding definitions of $P_m$ yield a term $\rho^{-1}P^{'}_m$ which,
by (\ref{eq:3.7}), contributes expressions $W_{jm}$ satisfying (\ref{eq:3.8}).
Now, although we have allowed the possibility that $W_{1m}$ may be zero,   there is 
no reason to suppose  that it is necessarily zero, and there is therefore nothing to be gained by departing from the simplest definition of $P_m$ based on (\ref{eq:1.5}) and (\ref{eq:1.6}).
\subsection{Other computational algorithms}
Here we indicate how our algorithm has a quite different purpose as compared to the algorithms of \cite{DC82} and \cite{D92}, and it is convenient to refer specifically to the latter.
In \cite{D92}, the differential system is
\begin{displaymath}
Y'(x)=x^{-1}B(x)Y(x)
\end{displaymath}
where $x$ can be a complex variable,
\begin{displaymath}
B=
\left (
\begin{array}{ccccc}
0 & 1 & &  &  \\
  & . & . & & \\
  &   & . & . & \\
  &   &  & 0 & 1 \\
-b_n & . & . &  -b_2 & -b_1 \\
\end{array}
\right )
\end{displaymath}
and each $b$ in the last row has a Laurent series
\beq
b(x)=({\rm const}.)x^c(1+a_1x^{-1}+...) \;\;\;(x \rightarrow \infty) \label{eq:6.1}
\enq
with rational $c$. With $\infty$ as an irregular singular point, the solutions of the corresponding $n-th$
order differential equation have the asymptotic form
\beq
f(x) \{ 1 + o(1) \}  \label{eq:6.2}
\enq
where the dominant term $f(x)$ comprises the usual logarithmic and exponential factors.
The algorithm developed by \cite{D92} determines $f(x)$ from a knowledge of $B$. The algorithm of ( Della Dora, Di Crescenzo and Tournier 1982)
also allows sub-dominant componants of $f(x)$ to be computed. 
In contrast, our algorithm is concerned with improvements to the $o(1)$ term in (\ref{eq:6.2}) and the construction of 
sub-dominant terms as explained above in section 1 and at the end of section 4.
It  is also the purpose of our paper to cover classes of coefficients which, unlike (\ref{eq:6.1}), contain   periodic factors as well as powers of $x$. This was the subject of section 5.
\par
We are grateful to the referees for raising the issues which are covered in this section.

\bibliographystyle{agsm}
\bibliography{bibliography}  
\end{document}